\documentclass[11pt]{amsart}
\usepackage{graphicx}
\usepackage{graphics}
\usepackage{amsmath}
\usepackage{amscd}
\usepackage{latexsym}
\usepackage[all]{xy}
\begin{document}
\textwidth 5.5in
\textheight 8.3in
\evensidemargin .75in
\oddsidemargin.75in

\newtheorem{lem}{Lemma}[section]
\newtheorem{conj}{Conjecture}[section]
\newtheorem{defi}{Definition}[section]
\newtheorem{thm}{Theorem}[section]
\newtheorem{cor}{Corollary}[section]
\newtheorem{lis}{List}[section]
\newtheorem{rmk}{Remark}[section]
\newtheorem{que}{Question}[section]
\newtheorem{prop}{Proposition}[section]
\newcommand{\p}[3]{\Phi_{p,#1}^{#2}(#3)}
\def\Z{\mathbb Z}
\def\R{\mathbb R}
\def\g{\overline{g}}
\def\odots{\reflectbox{\text{$\ddots$}}}
\newcommand{\tg}{\overline{g}}
\def\proof{{\bf Proof. }}
\def\ee{\epsilon_1'}
\def\ef{\epsilon_2'}
\title{A plug with infinite order and some exotic 4-manifolds}
\author{Motoo Tange}
\thanks{The author is partially supported by Grant-in-Aid for JSPS Fellows 21-1458}
\address{Research Institute for Mathematical Sciences, Kyoto University, Kyoto, 606-8502, Japan}
\email{tange@kurims.kyoto-u.ac.jp}
\date{\today}
\maketitle
\begin{abstract}
Every exotic pair in 4-dimension is obtained each other by twisting a {\it cork} or {\it plug} which are codimension $0$ submanifolds embedded in the 4-manifolds.
The twist was an involution on the boundary of the submanifold.
We define cork (or plug) with order $p\in {\Bbb N}\cup \{\infty\}$ and show
there exists a plug with infinite order.
Furthermore we show twisting $(P,\varphi^2)$ gives to enlargements of $P$ compact exotic manifolds with boundary.
\end{abstract}
\section{Introduction}
\subsection{Smooth structures}
Let $X$ be a smooth manifold.
If a smooth manifold $X'$ is homeomorphic but non-diffeomorphic to $X$,
then we say that $X$ and $X'$ are {\it exotic (pair)}.
Any exotic pair gives a different smooth structure on a topological manifold.
It is known that if $X$ has at least two smooth structures, then the dimension is greater than $3$.

Cork (or plug) is a pair of a submanifold with codimension 0 in a 4-manifold $X$ and an involution on the boundary.
They were defined in \cite{[A1]} and \cite{[GS]} and Akbulut-Yasui \cite{[AY1]}.
Twisting the cork (or plug) in $X$ by the involution, we can get an exotic pair.
Conversely it is known that any simply connected exotic pair can be obtained 
by the (contractible) cork twisting as proven in \cite{[AM],[CFHS],[M]}.

Recently many smooth structures have constructed by using cork and plug as in
\cite{[A2],[AY1],[AY2],[AY3],[AY4],[GS]}.
Since the main idea for the existence of cork and plug is due to the failure of the h-cobordism theorem in 4-dimension,
naturally the self-diffeomorphism of the boundary of cork and plug is an involution.
For example as appeared in \cite{[GS]} an exotic pair 
$E(2)\#\overline{{\Bbb C}P^2}$ and $\#^3{\Bbb C}P^2\#^{20}\overline{{\Bbb C}P^2}$ are
obtained by a cork twisting (an involution), that $C$ is a contractible 4-manifold having Mazur type and the involution is
by a symmetry of the framed link presenting $C$.

There exist {\it infinite} smooth structures as a character of 4-dimension.
Fintushel-Stern's {\it Knot surgery}, which is defined below, gives rise to mutually non-diffeomorphic manifolds.
Let $T\subset X$ be an embedded torus with trivial normal bundle.
Let $K$ be a knot in $S^3$.
The knot surgery is defined by
$$X_K=[X-\nu(T)]\cup_{\varphi}[(S^3-\nu(K))\times S^1],$$
where the definition of $\varphi$ is in \cite{[FS]}
and the notation $\nu$ stands for the tubular neighborhood of the submanifold.
The cut and paste notation will be defined in Definition~\ref{cp}.

One answer of the question is the result in \cite{[A2]}.
Our motivation is to construct ``a cork (or plug)" representing the infinite of exoticity.
We will relax the definition of cork and plug to accept infinite order.
In the next subsection we will define such a cork and plug and 
in next section we show the following.
\begin{thm}
\label{main}
There exists a plug $(P,\varphi)$ with infinite order.
$P$ is a simply connected, compact, Stein 4-manifold with $b_2=2$.
\end{thm}

The plug $P$ can be embedded in an elliptic fibration $X$ over $D^2$ with three vanishing cycles, which exactly two of them
are parallel.
In this case the plug twist $(P,\varphi)$ can change $X$ to knot surgery $X_K$,
where $K$ is any unknotting number $1$.
Furthermore $(P,\varphi^n)$ gives $X_{K^n}$, where $K^n$ is a knot obtained by $n$ times iteration
of the knotting operation from unknot to $K$.

\begin{figure}[ht]
\begin{center}
\includegraphics{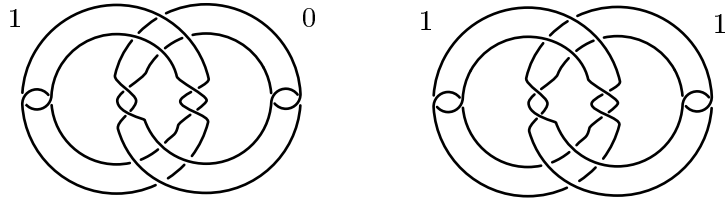}
\caption{$\tilde{Y}$ and $\tilde{Z}$.}
\label{boundary}
\end{center}
\end{figure}

The square $(P,\varphi^2)$ of the plug twist is a non-contractible cork with infinite order.
Using this we obtain the following.
\begin{thm}
\label{exoticex}
Let $\tilde{Y}$ be a 4-manifold presented by the left diagram in Figure~\ref{boundary}.
There exists a minimal (not having any $(-1)$-spheres), symplectic, simply-connected 4-manifold $Y_2$ with 
the same homeomorphism type as $Y=\tilde{Y}\#\overline{{\Bbb C}P^2}$.
In particular $Y$ and $Y_2$ are an exotic pair.

Let $\tilde{Z}$ be a 4-manifold presented by the right diagram in Figure~\ref{boundary}.
There exists a minimal (not having any $(-1)$-spheres), symplectic, simply-connected 4-manifold $Z_2$ with 
the same homeomorphism type as $Z=\tilde{Z}\#^2\overline{{\Bbb C}P^2}$.
In particular $Z$ and $Z_2$ are an exotic pair.

\end{thm}

Finally we define a notation on cut and paste of 4-manifolds, used in the paper.
\begin{defi}[Cut and Paste]
\label{cp}
We define a notation of cut and paste.
Let $X$ be a 4-manifold and $Z\subset X$ a compact, codimension-$0$ submanifold in $X$.
Let $Y$ be a 4-manifold with the same boundary $M$ as $\partial Z$.
We fix identification $M=\partial Z=\partial Y$.
Let $\varphi$ be a self-diffeomorphism of $M$.
The manifold obtained by identifying $\partial W$ and $\partial Y$ through $\partial Z=M\overset{\varphi}{\to}M=\partial Y$ 
is
presented as follows:
$$[X-Z]\cup_{\varphi} Y.$$
If $Z=Y$, we use the obvious map as the identification.
\end{defi}
\subsection{Cork and Plug with order $p$.}
We define notions of cork and plug with order $p$.
See \cite{[AY1]} for the original cork and plug.
\begin{defi}[Cork with order $p$]
Let $(C,\varphi)$ be a pair of a compact, contractible, Stein 4-manifold $C$ with boundary and a diffeomorphism $\varphi:\partial C\to \partial C$.
The pair $(C,\varphi)$ is called a cork with order $p$ if $C$ satisfies the following properties:
the order of $\varphi$ is $p$($\ge 2$), $\varphi$ can extend to a self-homeomorphism of $C$ but
$\varphi^q\ (1\le {}^{\forall}q< p)$ cannot be extended to any self-diffeomorphism of $C$.

In the case of $\varphi^q\neq$id for any natural number $q$ we call $(C,\varphi)$ a cork with infinite order.

Let $(C,\varphi)$ be a cork with order $p$ and $X$ a 4-manifold
containing $C$.
A cork twist of $X$ is the set of manifolds $[X-C]\cup_{\varphi^q} C$\ \ ($1\le q\le p$).
If the cork twist of $X$ gives mutually different $p$ smooth structures,
the cork $(C,\varphi)$ is called a cork of $X$ with order $p$.
\end{defi}
When we treat a non-contractible submanifold as $C$ as appeared in \cite{[AY1]}, then
we call it {\it a generalized cork with order $p$}.

\begin{defi}[Plug with order $p$]
Let $(P,\varphi)$ be a pair of a compact Stein 4-manifold $P$ with boundary and a diffeomorphism $\varphi:\partial P\to \partial P$.
The pair $(P,\varphi)$ is called a plug with order $p(\ge 2)$ if $P$ satisfies the following properties:
the order of $\varphi$ is $p$($\ge 2$), $\varphi$ cannot be extended to any self-homeomorphism of $P$,
there exists a 4-manifold $X$ containing $P$ such that $[X-P]\cup_{\varphi^q}P$ ($0\le {}^{\forall}q< p)$ are
$p$ mutually non-diffeomorphic manifolds.

In the case of $\varphi^q\neq$id for any natural number $q$ we call $(P,\varphi)$ a plug with infinite order.

Let $(P,\varphi)$ be a plug with order $p$ and $X$ a 4-manifold
containing $P$.
A plug twist of $X$ is the set of manifolds $[X-P]\cup_{\varphi^q} P$\ \ ($0\le  q< p$).
If the plug twist of $X$ gives mutually different $p$ smooth structures,
the plug $(P,\varphi)$ is called a plug of $X$ with order $p$.
\end{defi}
Any cork (or plug) with order $2$ means the original cork (or plug).

\subsection{Acknowledgements}
I thank Professor Kouichi Yasui and Yuichi Yamada for giving me some useful and suggestion and advice to study
plug with infinite order and also thank Professor Mikio Furuta for a useful comment in the meeting "Four-Dimensional Topology in 2011".

\section{A plug $(P,\varphi)$ with infinite order.}
\subsection{The diffeomorphism type of $P$.}
\begin{figure}[ht]
\begin{minipage}{.495\hsize}
\begin{center}
\includegraphics{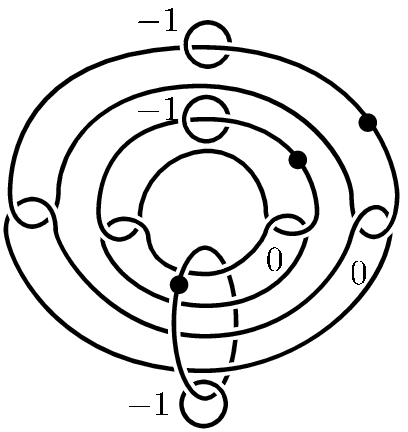}
\caption{A handle decomposition of $P$.}
\label{P}
\end{center}
\end{minipage}
\begin{minipage}{.495\hsize}
\begin{center}
\includegraphics{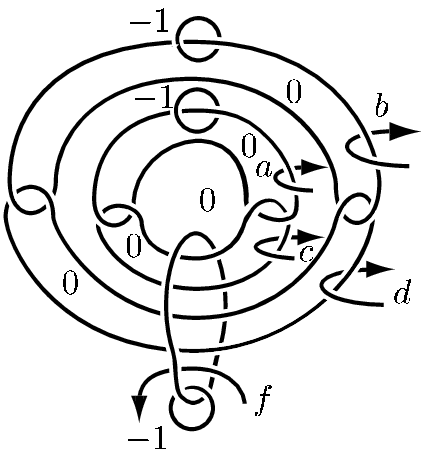}
\caption{The generators in $\pi_1(\partial P)$.}
\label{gene}
\end{center}
\end{minipage}
\end{figure}

We define a plug $P$ mainly used through the paper.
\begin{defi}
\label{defin}
We define a compact 4-manifold $P$ to be a manifold admitting the handle decomposition in Figure~\ref{P}.
\end{defi}
On the other hand the manifold $P$ is obtained as $V_2\times S^1$ with three $-1$-framed 2-handles,
where $V_2$ is the genus 2 handlebody.

Thus $P$ is simply connected and has $H_2(P)\cong {\Bbb Z}^2$.
The fundamental group of the boundary $\partial P$ is
$$\pi_1(\partial P)=\langle a,b,c,d,f|[f,b^{-1}],[f,a^{-1}],a[f,c^{-1}],b[f,d^{-1}],f^{-1}YX\rangle,$$
where we use generators $a,b,c,d$ and $f$ in Figure~\ref{gene}.
Elements $X,Y$ are $X=[b^{-1},d]$ and $Y=[a^{-1},c]$.
Therefore we get $H_1(\partial P)\cong{\Bbb Z}^2$.

\subsection{A diffeomorphism $\varphi$ on $\partial P$.}
\begin{figure}[htbp]
\begin{center}
\includegraphics{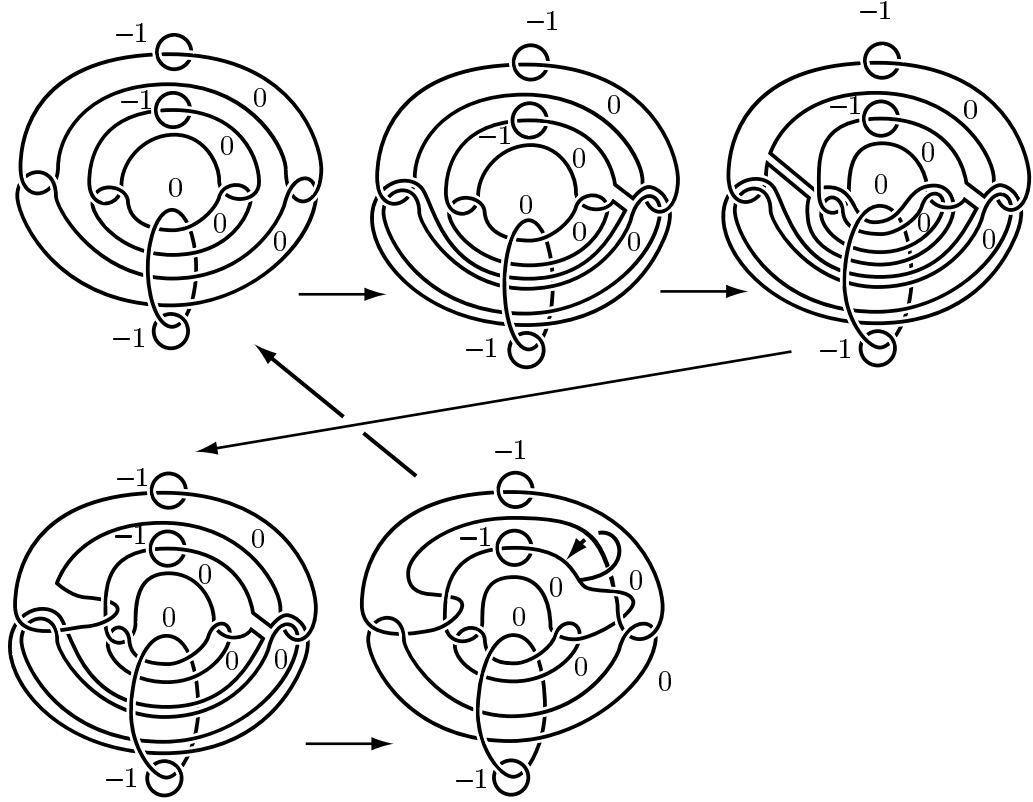}
\caption{A diffeomorphism $\varphi$.}
\label{Pdiff}
\end{center}
\end{figure}

We define a diffeomorphism on $\partial P$.
Moving the diagram of $\partial P$ in accordance with the process of Figure~\ref{Pdiff}
obtains a diffeomorphism of $\partial P$ is denoted by $\varphi$.

From the definition of $\varphi$ the images of the generators $a,b,c,d,f$ are as follows:
\begin{equation}\label{indu}
(Xb)^{-1}aXb,\ b,\ (Xb)^{-1}a(Xb)a^{-1}c(Xb),\ d(Xd)^{-1}a(Xd),\ f.
\end{equation}

\begin{lem}
The gluing map $\varphi$ acts on $H_1(\partial P)$ and $H_2(\partial P)$ trivially.
\end{lem}
\proof
The images in (\ref{indu}) on the abelianization of $\pi_1(\partial P)$ are trivial.
Then the abelianization of $\varphi_\ast$ is trivial.
Thus $\varphi_\ast$ acts on $H_2(\partial P)$ trivially through the Poincar\'e duality.

\subsection{$P$ is a Stein manifold.}
\begin{figure}[ht]
\begin{center}
\includegraphics{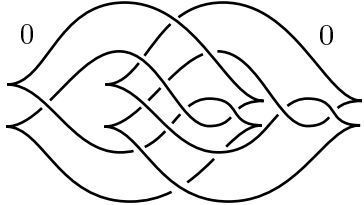}
\caption{A Stein structure over $P$.}
\label{PStein}
\end{center}
\end{figure}

To show that $P$ is a plug, $P$ must be a Stein manifold  (admit a Stein structure).
Admitting 4-dimensional Stein structure is due to a description by a specific Legendrian surgery diagram on $\#^nS^2\times S^1$.
The manifold is constructed by attaching 2-handles on $\natural^nD^3\times S^1$ along the framed Legendrian link,
where $\natural$ stands for the boundary sum.
The condition of the framed Legendrian link is that each framing of attaching 2-handles are $\text{tb}(K)-1$,
where $\text{tb}(K)$ is the Thurston-Bennequin invariant of the Legendrian knot $K$
(see \cite{[GS]} for the explanation).
\begin{prop}
\label{stein}
$P$ admits a Stein structure.
\end{prop}
\proof
$P$ admits handle decomposition as in Figure~\ref{PStein}.
Each Thurston-Bennequin invariant of the components is $1$.
Thus all framings of the components satisfy $\text{tb}(K)-1$.
Therefore $P$ admits a Stein structure.\qed
\subsection{Infinite exotic manifolds from $(P,\varphi)$-twist.}
\begin{figure}[htbp]
\begin{center}
\includegraphics{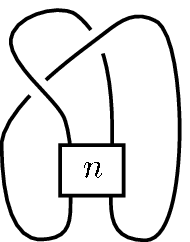}
\caption{$K_n$, where \fbox{$n$} is a full $n$-twist.}
\label{twistknot}
\end{center}
\end{figure}

\begin{figure}[htbp]
\begin{center}
\includegraphics{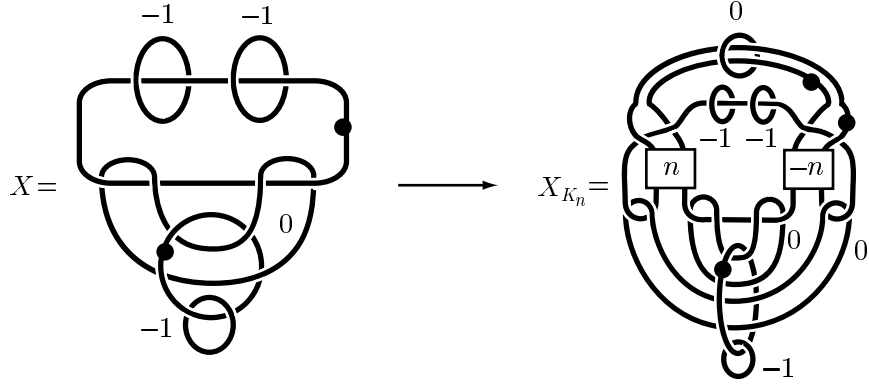}
\caption{The knot surgery of $X$.}
\label{X}
\end{center}
\end{figure}

This subsection is essential for the pair $(P,\varphi)$ to be a plug with infinite order.
For any embedding $P\subset W$ the performance $[W-P]\cup_\varphi P$ by using $(P,\varphi)$ is called $(P,\varphi)$-twist.
\begin{prop}
\label{infinity}
The pair $(P,\varphi)$ produces infinitely many exotic manifolds.
\end{prop}
\proof
Let $K_n$ be a twist knot as in Figure~\ref{twistknot}.
Let $X$ be a manifold as in the left diagram in Figure~\ref{X}.
Performing knot surgery on $X$, we get the right diagram in Figure~\ref{X}.
Sliding one of the two $-1$-framed 2-handles over $0$-framed 2-handle in the diagram, and removing the top $0$-framed 2-handle, we get the
handle diagram of $P$ in $X$.

\begin{figure}[htbp]
\begin{center}
\includegraphics{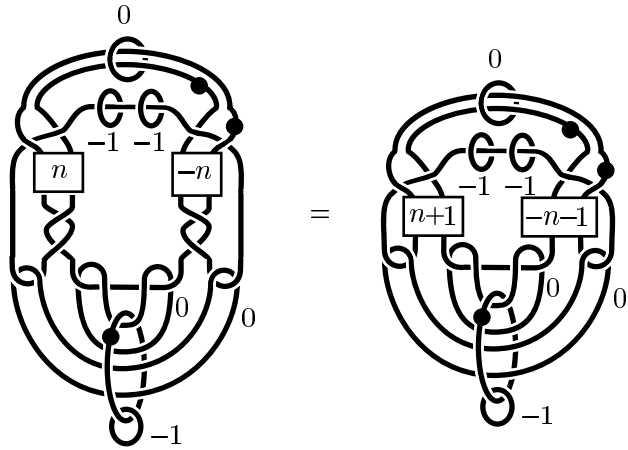}
\caption{The $(P,\varphi)$ twist $[X_{K_n}-P]\cup_\varphi P$ is diffeomorphic to $X_{K_{n+1}}$.}
\label{XX}
\end{center}
\end{figure}

Here we perform the $(P,\varphi)$-twist $[X_{K_n}-P]\cup_\varphi P$.
Keeping track of the diffeomorphisms in Figure~\ref{Pdiff} for $X_{K_n}$,
consequently from Figure~\ref{XX} we get the following:
$$[X_{K_n}-P]\cup_{\varphi}P=X_{K_{n+1}}.$$
Namely $(P,\varphi)$-twist of $X_{K_n}$ gives $X_{K_{n+1}}$.

Here we embed $X$ in $E(2)$ (the K3-surface) as an elliptic fibration having three vanishing cycles
which exactly two of them are parallel.

We trivialize the tubular neighborhood $E(2)\supset \nu(T)=D^2\times T^2$ of the general fiber $T$,
where the vanishing cycles of $E(2)$ generate homology classes of $\{\text{pt}\}\times T$ and the direction 
$\partial D^2\times \{\text{pt}\}$ corresponds to a section of $E(2)$.
This section is mapped to the longitude of $K_n$ by the gluing map of knot surgery.

\begin{figure}[htbp]
\begin{center}
\includegraphics{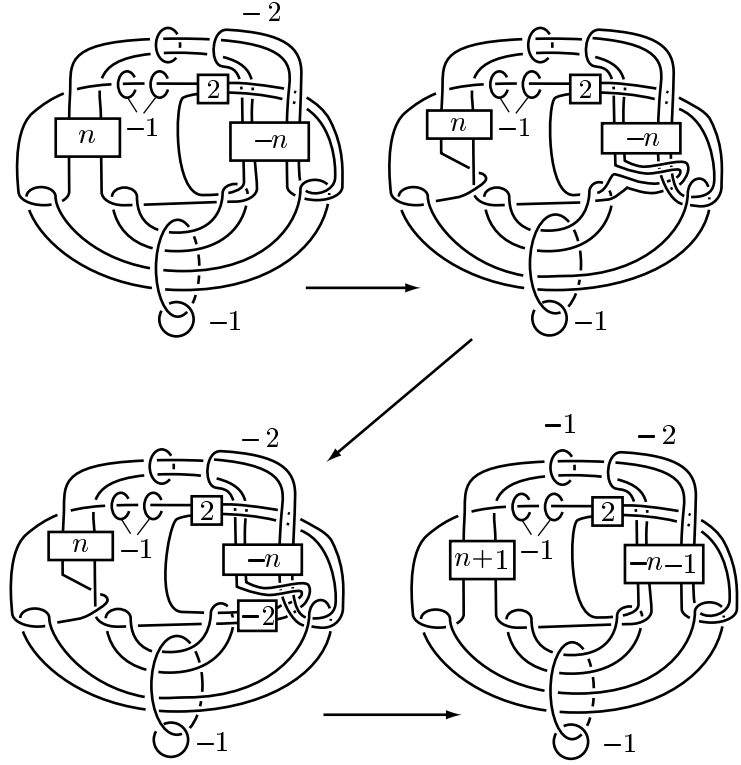}
\caption{The two longitudes preserves by the operation.}
\label{section}
\end{center}
\end{figure}

The obvious map $i:\partial X_{K_n}\to \partial(X_{K_n}-P)\cup_\varphi P=\partial X_{K_{n+1}}$ preserves
the two directions corresponding to the vanishing cycles (see Figure~\ref{Pdiff}).
Furthermore $i$ maps the direction of the section in $\partial X_{K_n}$  (namely longitude of $K_n$) to
the direction of the section in $\partial X_{K_{n+1}}$ (see Figure~\ref{section}, where 
the framings of unlabeled links are $0$).
Thus in this case the replacement of $P$ by $\varphi$ of $E(2)$ gives $E(2)_{K_{n+1}}$.
$$[E(2)_{K_{n}}-P]\cup_{\varphi}P=E(2)_{K_{n+1}}.$$
In general for any 4-manifold $E$ containing $X$ the replacement of $P$ by $\varphi$
gives some knot surgery $E_K$.
$K$ is determined depending on $P\hookrightarrow X$.

By the Seiberg-Witten formula on knot surgery in \cite{[FS]}
$$SW_{E(2)_{K_n}}=nt-2n-1+nt^{-1},$$
thus $E(2)_{K_n}$ are infinitely many mutually non-diffeomorphic manifolds by $(P,\varphi^i)$-twist.
As a result $(P,\varphi)$-twist can give rise to infinite exotic pairs.
\qed
\medskip\\
Namely the map $\varphi$ has infinite order.
\subsection{Extendability as a homeomorphism.}
\begin{figure}[htbp]
\begin{center}
\includegraphics{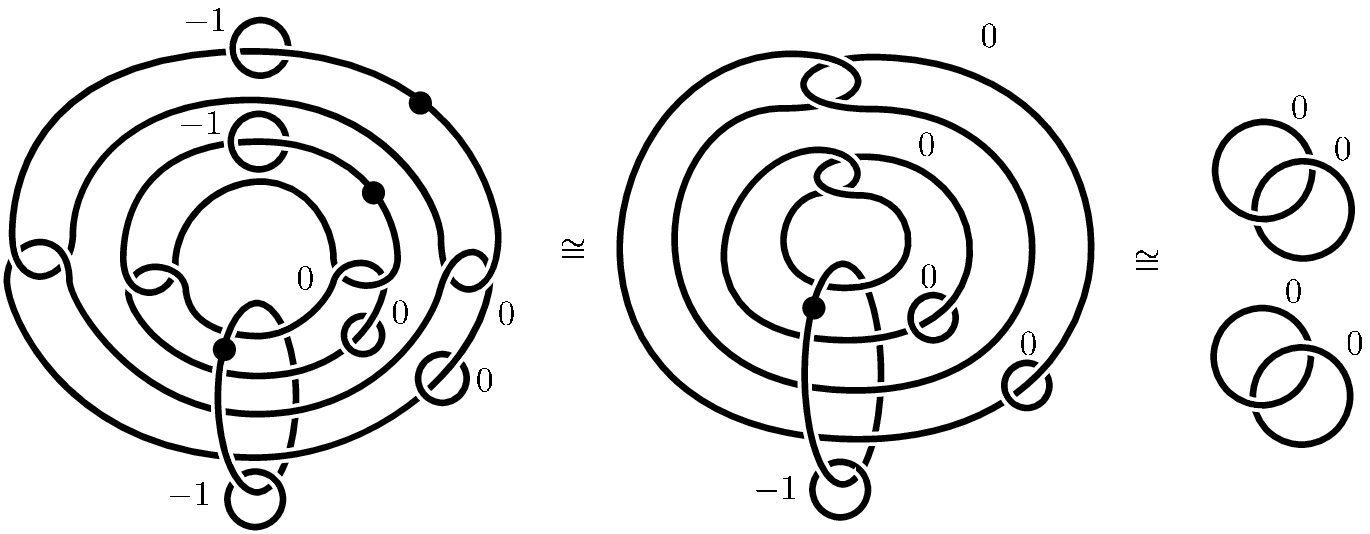}
\caption{$D(P)=P\cup_{\text{id}}\overline{P}$.}
\label{double}
\end{center}
\end{figure}

\begin{figure}[htbp]
\begin{center}
\includegraphics{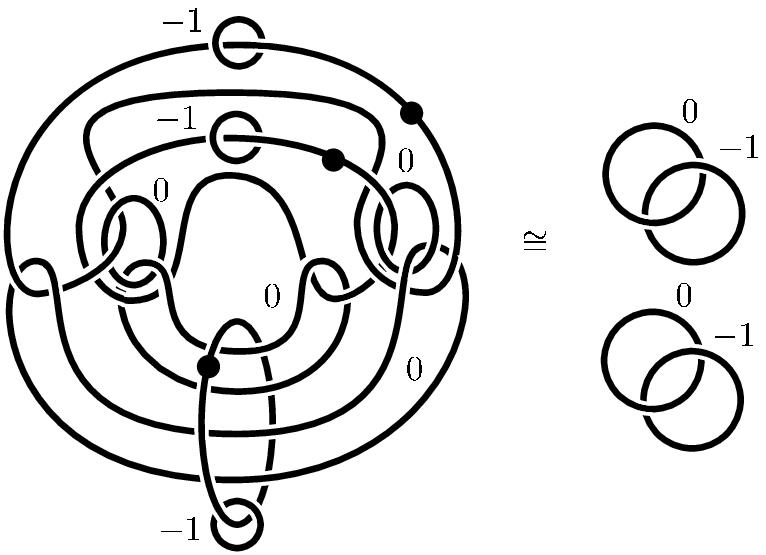}
\caption{$D(P)_{\varphi}=P\cup_{\varphi}\overline{P}$.}
\label{double2}
\end{center}
\end{figure}

Here we show that $\varphi$ cannot be extended to a homeomorphism of $P$.
\begin{prop}
\label{extend}
Let $n$ be any integer.
The diffeomorphism $\varphi^n:\partial P\to \partial P$
$$\left\{
\begin{array}{ll}
\text{cannot be extended to $P\to P$ as any homeomorphism}&n\text{ odd}\\
\text{can be extended to $P\to P$ as a homeomorphism}&n\text{ even.}\\
\end{array}
\right.$$
In particular $P\cup_{\text{id}}\overline{P}$ and $P\cup_{\varphi}\overline{P}$ are
diffeomorphic to $\#^2S^2\times S^2$ and $\#^2({\Bbb C}P^2\#\overline{{\Bbb C}P^2})$ respectively.
\end{prop}
\proof
We take the double $D(P):=P\cup_{\text{id}} \overline{P}$ by using the identity map.
This manifold is diffeomorphic to $\#^2S^2\times S^2$ due to Figure~\ref{double}.
On the other hand the manifold $D(P)_\varphi:=P\cup_{\varphi}\overline{P}$ glued by use of $\varphi$
is diffeomorphic to $\#^2({\Bbb C}P^2\#\overline{{\Bbb C}P^2})$ as in Figure~\ref{double2}.
Since $\#^2S^2\times S^2$ and $\#^2({\Bbb C}P^2\#\overline{{\Bbb C}P^2})$ are not homeomorphic,
then $\varphi$ cannot be extended to any homeomorphism of $P$.


\begin{figure}[htbp]
\begin{center}
\includegraphics{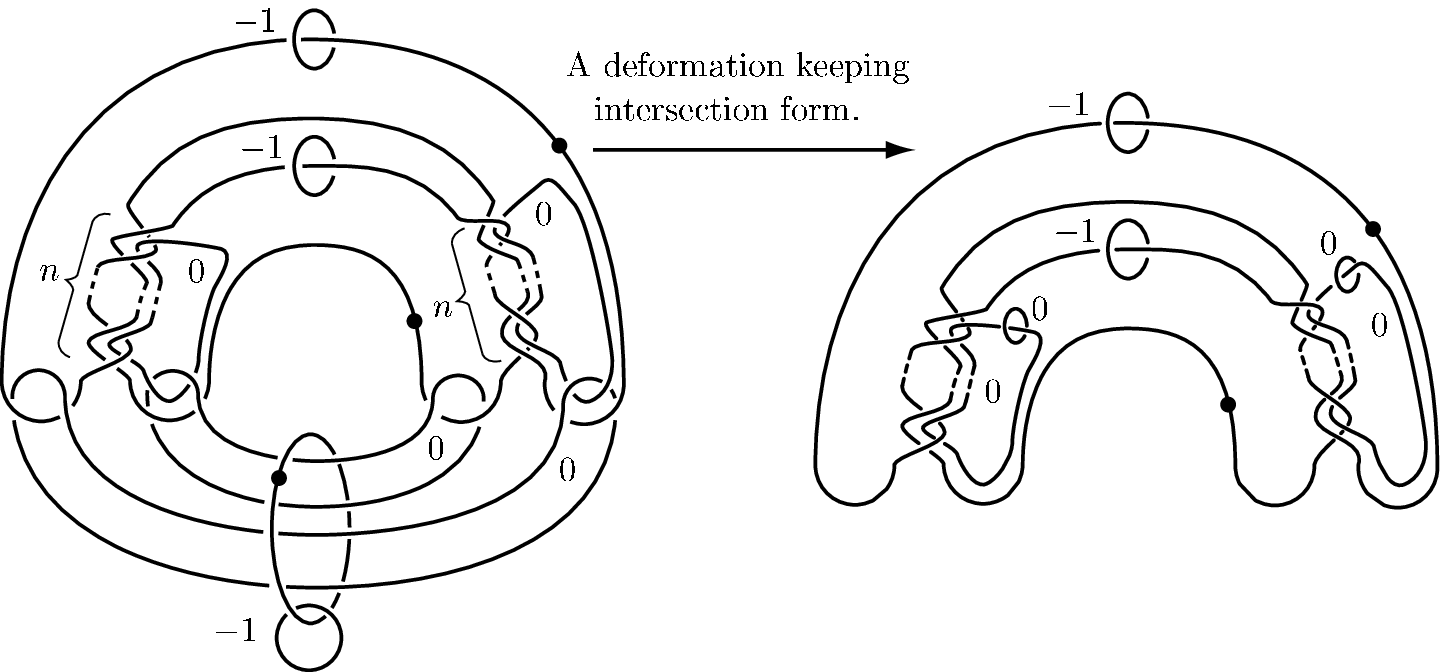}
\caption{$D(P)_{\varphi^n}\to D_n$ (homeomorphism).}
\label{doublen}
\end{center}
\end{figure}

\begin{figure}[htbp]
\begin{center}
\includegraphics{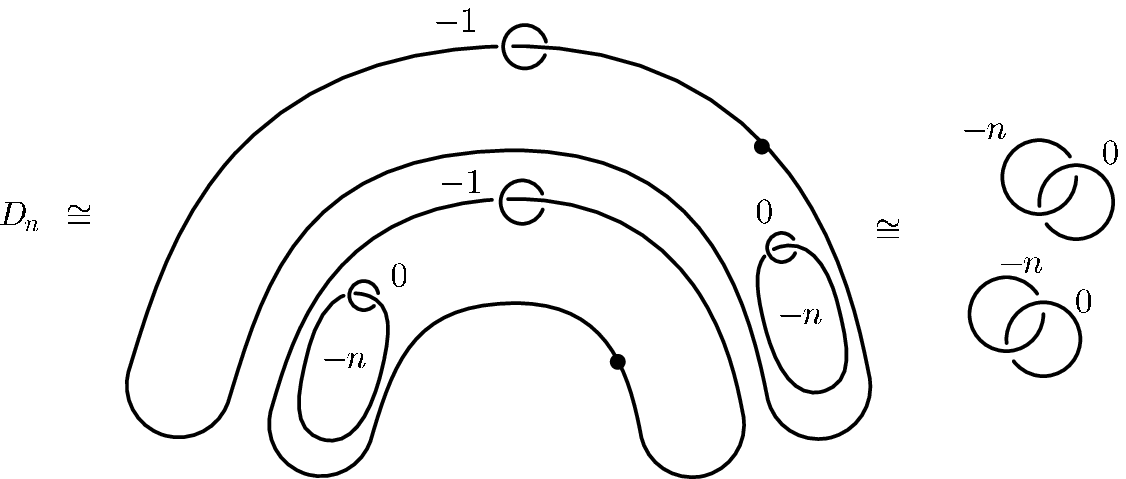}
\caption{Diffeomorphism types of $D_n$.}
\label{mwari2}
\end{center}
\end{figure}

In general $D(P)_{\varphi^n}=P\cup_{\varphi^n}\overline{P}$ is diffeomorphic to the left in Figure~\ref{doublen}.
Replacing $-1$-framed 2-handle at the bottom in this picture with $0$-framed 2-handle,
we get a manifold $D_n$ with the same intersection form as $D(P)_{\varphi^n}$.
From Freedman's result $D(P)_{\varphi^n}$ and $D_n$ are homeomorphic.
\begin{figure}[htbp]
\begin{center}
\includegraphics{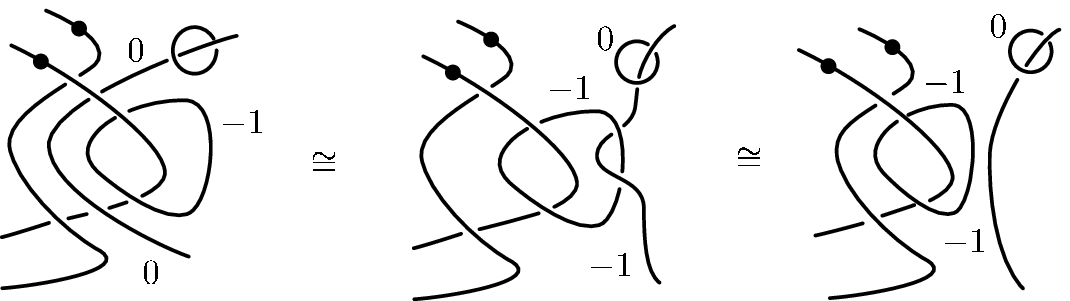}
\caption{A diffeomorphism.}
\label{mwari}
\end{center}
\end{figure}
The diagram of $D_n$ in Figure~\ref{doublen} can be simplified by iterating the local process in Figure~\ref{mwari}.
As in Figure~\ref{mwari2}, $D_n$ is diffeomorphic to 
$$\begin{cases}
\#^2({\Bbb C}P^2\#\overline{{\Bbb C}P^2})&n\text{ odd}\\
\#^2S^2\times S^2&n\text{ even.}
\end{cases}$$
Hence $D_{\varphi^n}(P)$ is homeomorphic to $\#^2S^2\times S^2$ or 
$\#^2({\Bbb C}P^2\#\overline{{\Bbb C}P^2})$ respectively 
namely it is spin (or non-spin) if $n\equiv 0(2)$ (or $n\equiv 1(2)$).

According to (0.8) Proposition (iii) in Boyer's paper \cite{[B]},
$\varphi^n$ can be extended to a homeomorphism of $P$ if $n$ is even,
and $\varphi^n$ cannot be extended to any homeomorphism of $P$ if $n$ is odd.
\qed
\medskip\\
As a corollary we get the following corollary.
\begin{cor}
\label{cork1}
For non-zero integer $m$,
$(P,\varphi^{2m})$ are infinite generalized corks with infinite order.
\end{cor}
\proof
Assertion in Proposition~\ref{extend} means $\varphi^{2m}$
can be extended as a homeomorphism of $P$.
Therefore $(P,\varphi^{2m})$ are infinite many corks with infinite order.
\qed
\medskip

{\bf Proof of Theorem~\ref{main}.}\\
Definition~\ref{defin}, Proposition~\ref{stein},~\ref{infinity}, and ~\ref{extend} mean Theorem~\ref{main}.\qed
\medskip

From the proof of Proposition~\ref{infinity}
the plug twist by $(P,\varphi)$ means ``crossing changing operation" of smooth structures through knot surgery.
On the other hand any twist knot $K_n$ is an unknotting number $1$ knot and there exist two embeddings $\iota_1:P\hookrightarrow E(2)$
and $\iota_2:P\hookrightarrow E(2)$ such that the plug twists give rise to $E(2)_{K_n}$ and $E(2)_{K_{n+1}}$ as
the diagram below.
$$
\xymatrix{
E(2) \ar@{->}[r]^{(\iota_1,P,\varphi)}\ar@{->}[dr]_{(\iota_2,P,\varphi)}  &  E(2)_{K_n}\ar@{->}[d]^{(P,\varphi)}\\
 & E(2)_{K_{n+1}}\\
}
$$
This means $\iota_1$ and $\iota_2$ are different embeddings of $P$ in $E(2)$.

We raise the several questions.
\begin{que}
Are there exist any plug $(Q,\psi)$ with infinite order which $\psi^n$ for any integer $n$ cannot be extended to a homeomorphism of $P$?

Or if $(Q,\psi)$ is a plug with infinite order, then is $(Q,\psi^2)$ a generalized cork with infinite order?
\end{que}

\begin{que}
Are there exist any plug or cork $(Q,\psi)$ with finite order $p$ $(3\le p<\infty)$?
\end{que}
\section{Some exotic manifolds.}
We consider enlargements of $P$ attaching some $-1$-framed 2-handles.
\subsection{A manifold exotic to $Y$.}
\begin{figure}[htbp]
\begin{center}
\includegraphics{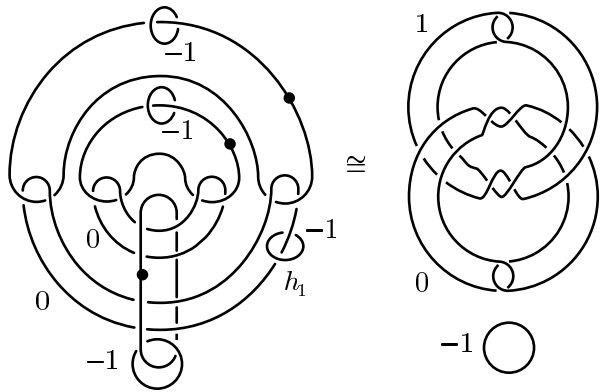}
\caption{$Y:=P\cup h_1=\tilde{Y}\#\overline{{\Bbb C}P^2}$.}
\label{PH}
\end{center}
\end{figure}

Attaching $-1$-framed 2-handle over one meridian of the two $0$-framed 2-handles of $P$
as in Figure~\ref{PH} we define the resulting manifold to be $Y$.
Thus $Y:=\tilde{Y}\#\overline{{\Bbb C}P^2}$, where
$\tilde{Y}$ is the 4-manifold attached along a satellite link as in right picture in Figure~\ref{PH}.

$\tilde{Y}$ is a simply connected 4-manifold with $b_2=2$ and
the boundary $\partial Y$ is a 3-manifold with $H_1\cong {\Bbb Z}$.

Let $W$ be a 4-manifold and $T\subset W$ an embedded torus with the trivial neighborhood.
Let $L=K_1\cup\cdots\cup K_n$ an $n$-component link.
For $n$ copies $(W_i,T_i)$ of $(W,T)$ we define {\it link operation} $W_L$ as follows:
$$W_L=[(S^3-\nu(L))\times S^1]\cup_{\varphi}\coprod_{i=1}^n[W_i-\nu(T_i)],$$
where the $n$ attaching maps are the common $\varphi$ and the definition is in \cite{[FS]}.
In this section we show the following.

\begin{figure}[htbp]
\begin{minipage}{.42\textwidth}
\begin{center}
\includegraphics{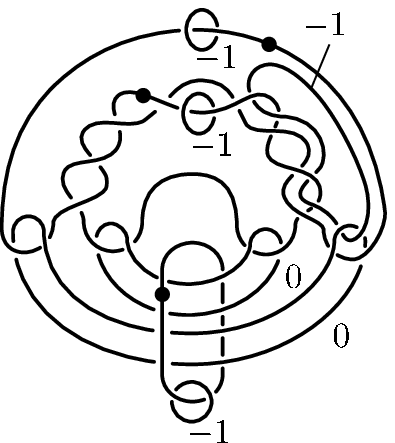}
\caption{$Y_2$.}
\label{FF}
\end{center}
\end{minipage}
\begin{minipage}{.42\textwidth}
\begin{center}
\includegraphics{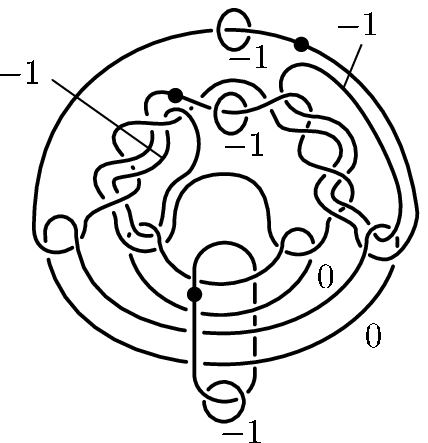}
\caption{$Z_2$.}
\label{F}
\end{center}
\end{minipage}
\end{figure}

\begin{thm}
\label{goon}
$Y$ admits at least two smooth structures $\{Y,Y_2\}$.
$Y_2$ admits symplectic structure and minimal.
\end{thm}
\proof
Twisting $\tilde{P}$ using the generalized cork $(P,\varphi^2)$, we get Figure~\ref{F} (called $Y_2$ here).
The manifold $Y_2$ can be embedded in a link operation $E(1)_{L_2}$, where $E(1)$ is the elliptic fibration with $12$ nodal singularities over $S^2$
and contains a general fiber as an embedded torus with trivial normal bundle.
$L_n$ is the $(2,2n)$-torus link.
Since the boundaries are homeomorphic by an easy handle calculation and 
the cork can be extended to $P$ as a homeomorphism, the two manifolds $Y$ and $Y_2$ are homeomorphic.

Since the Seiberg-Witten invariant of $E(1)_{L_2}$ is computed as the multivaluable Alexander polynomial
of $L_2$, it is
$$t_1t_2+t_1^{-1}t_2^{-1},$$
where each $t_i$ is the Poincar\'e dual $PD([T_i])$ of a general fiber $T_i$ of the two $E(1)$.
Namely the basic classes of $E(1)_{L_2}$ are $\{\pm PD([T_1]+[T_2])\}$.
If $Y$ and $Y_2$ are diffeomorphic, then there exists an embedding $Y\hookrightarrow E(1)_{L_2}$.
This means $E(1)_{L_2}=Y'\#\overline{{\Bbb C}P^2}$.
The Seiberg-Witten basic classes $\mathcal{B}_{E(1)_{L_2}}$ of $E(1)_{L_2}$ have to be of form $\mathcal{B}_{E(1)_{L_2}}=\{K\pm PD([E_1])|K\in \mathcal{B}_{Y'}\}$, where $E_1$ are 
the exceptional sphere.
Thus $[E_1]=\pm ([T_1]+[T_2])$ holds.
From the square $-1=[E_1]^2=([T_1]+[T_2])^2=0$ this is contradiction.
Therefore $Y$ and $Y_2$ are non-diffeomorphic and $Y_2$ does not have any $(-1)$-sphere.
Thus $Y$ and $Y_2$ are exotic manifolds.

Furthermore $L_2$ is a fibered link, thus $E(1)_{L_2}$ is a symplectic manifold.
In particular $Y_2$ is a minimal symplectic manifold.
\qed
\medskip

Applying the same argument to the case of the generalized cork twist $(P,\varphi^{2n})$ of $Y$,
we obtain $Y_{2n}\hookrightarrow E(1)_{L_{2n}}$.
$$SW_{E(1)_{L_{2n}}}=\Delta_{L_{2n}}(t_1,t_2)=(t_1t_2)^{2n-1}+(t_1t_2)^{2n-3}+\cdots+(t_1t_2)^{-2n+1}$$
As a result $Y$ and $Y_{2n}$ are non-diffeomorphic.
Since each $Y_{2n}$ does not contain ($-1$)-sphere, $Y_{2n}$ is a minimal symplectic 4-manifold.
However whether these manifolds $Y_{2n}\ (n\ge 1)$ are mutually non-diffeomorphic or not is not known.
\medskip

\subsection{A manifold exotic to $Z$.}
Next we define $Z=P\cup h_1\cup h_2=\tilde{Z}\#^2\overline{{\Bbb C}P^2}$,
where the attaching circles of $h_1$ and $h_2$ are the meridians of $0$-framed 2-handles in Figure~\ref{PH}.
The framings of $h_1$ and $h_2$ are $-1$.
$\tilde{Z}$ is a simply connected manifold obtained by attaching two 2-handles over the same link as $\tilde{Z}$ with framing $\{1,1\}$.
$Z$ has $b_2=4$ and $b_2^+=2$.
The boundary $\partial Z$ is a homology sphere (a torus sum of two copies of the trefoil complement).

In the same way we assert the following.
\begin{thm}
\label{Z}
$Z$ admits at least two smooth structures $Z, Z_2$.
$Z_2$ admits symplectic structure and minimal.
\end{thm}
\proof
The manifold $Z_2$ which is obtained by the generalized cork twist $(P,\varphi^{2})$ can be 
embedded in $E(1)_{L_2}$.
Therefore by Freedman's result, $Z$ and $Z_2$ are homeomorphic.
In particular they have the same boundary.

If $Z_2$ is diffeomorphic to $Z$, the same argument of the Seiberg-Witten invariant as the
previous one leads to contradiction.
Therefore $Z_2$ is non-diffeomorphic to $Z$.
\qed

{\bf Proof of Theorem~\ref{exoticex}.}
Theorem~\ref{goon} and Theorem~\ref{Z} imply Theorem~\ref{exoticex}.\qed

\begin{conj}
$\{Y_{2m}\}$ and $\{Z_{2m}\}$ are infinitely many mutually exotic manifolds to $Y$ and $Z$ respectively.
\end{conj}


\end{document}